\documentclass[a4paper]{article}

\usepackage{MyLaTeX,MyArticle}
\RequirePackage{pgf,pgffor}
\usepackage{tikz}
\tikzstyle{every node}=[circle, fill=black!0, inner sep=0pt, minimum width=4pt]

\newcommand{\Deg}{\mathrm{Deg}}
\newcommand{\aA}{\mathrm{aA}}
\newcommand{\Arg}{\mathrm{Arg}}
\title{The Brauer Trees of Non-Crystallographic Groups of Lie Type}
\author{David A.~Craven, University of Birmingham}
\date{October 2012}

\begin{document}
\maketitle

\noindent In this article we determine the Brauer trees of the unipotent blocks with cyclic defect group in the `groups' $I_2(n,q)$, $H_3(q)$ and $H_4(q)$. The degrees of the unipotent characters of these objects were given by Lusztig, and using the general theory of perverse equivalences we can reconstruct the Brauer trees that would be consistent with Deligne--Lusztig theory and the geometric version of Brou\'e's conjecture. We construct the trees using standard arguments whenever possible, and check that the Brauer trees predicted by Brou\'e's conjecture are consistent with both the mathematics and philosophy of blocks with cyclic defect groups.

\medskip


\section{Introduction}

The representation theory of the finite groups of Lie type $G=G(q)$ is a vast area, much of it controlled in various ways by the associated Weyl group $W$. If $W$ is replaced by a finite Coxeter group that is \emph{not} a Weyl group, i.e., $H_3$, $H_4$ or $I_2(n)$, then there is of course no associated finite group of Lie type, $H_3(q)$ for example.

This has not stopped various authors from trying to develop a representation theory of the `group'. In 1993, Lusztig \cite{lusztig1993} gave a collection of unipotent degrees (polynomials in $q$), which should mimic the degrees of the unipotent characters of groups of Lie type. In 1997, M\"uller \cite{muller1997} determined the decomposition matrices of the corresponding Hecke algebras, which form submatrices of the decomposition matrices of the unipotent blocks. More recently, the theory of \emph{spetses} has emerged, by Brou\'e, Malle and Michel \cite{bmm1999} \cite{bmm2012un}, which attempt to construct an object associated even to each complex reflection group, which has in particular a set of unipotent degrees.

In a parallel development to the theory of spetses, \emph{perverse equivalences} \cite{chuangrouquierun} have been defined, and their use in the theory of groups of Lie type has started to be studied \cite{cravenrouquier2010un}. For blocks with cyclic defect groups, the structure of the cyclotomic Hecke algebra of a unipotent block was proved by the author in \cite{craven2012un} to completely control the structure of the Brauer tree of the block, in the sense that the tree can be reconstructed from the cyclotomic Hecke algebra. (There are currently three potential exceptions to this statement, two unipotent blocks of $E_8(q)$ and one of ${}^2\!F_4(q^2)$, and are exactly the unipotent blocks for which the labelling of the Brauer tree is not known.)

The non-crystallographic Coxeter groups have cyclotomic Hecke algebras attached to their `unipotent blocks' as well, so the methods of \cite{craven2012un} can be applied to these groups, and candidate Brauer trees for the blocks can be constructed. (This has led to the construction of all Brauer trees for unipotent blocks for groups of Lie type and all primes in \cite{cdr2012un}.) Although we can of course not \emph{prove} that these are the Brauer trees for the blocks, since the blocks do not actually exist, we can prove in many cases that the Brauer trees given here are the only possibility, or offer significant evidence that they are correct, for example, analogous results holding for groups of Lie type, and so on.

We echo the sentiments of Lusztig in \cite{lusztig1993}, in that the results in this paper are ``not theorems in the accepted sense", as they are still developing the theory of an object that does not exist. Nevertheless, the results here prove that the combinatorial version of Brou\'e's conjecture -- i.e., the derived equivalence implied by the geometric version of Brou\'e's conjecture suggested by Brou\'e and Malle \cite{brouemalle1993} -- as stated in \cite{craven2012un} but see also Section \ref{sec:combbroue}, carries through perfectly to the cyclic-defect blocks for the non-crystallographic case.

The current state of the combinatorial Brou\'e conjecture is such that, assuming it, the decomposition matrices of all unipotent blocks can be computed, but it is not possible to determine the Morita equivalence class of the unipotent blocks for non-cyclic defect, as it requires knowledge of the Green correspondents of the simple modules in a block. Since the Green correspondents are not defined in the non-crystallographic case, we cannot complete this task. However, if in a future version of the combinatorial Brou\'e conjecture this dependence is removed, as is likely, then we should be able to construct the representation theory of these `groups' for all primes.

\section{Notation and setup}

Write $\zeta_n$ for the primitive $n$th root of unity $\e^{2\pi\I/n}$. Let $q$ be a power of a prime $p$, and let $\ell\neq p$ be another prime. If $G=G(q)$ is a group of Lie type and $\ell\mid |G(q)|$, then there exists an integer $d$ such that $\ell\mid\Phi_d(q)$, since the order of $G(q)$ is a product of cyclotomic polynomials  $\Phi_d(q)$ and a power of $q$. The `generic' theory of $G(q)$ is when there is a single $d$ such that $\ell\mid\Phi_d(q)$, in which case the Sylow $\ell$-subgroup of $G$ is abelian. (One is able to do slightly better than this situation if one considers non-principal blocks, but this situation suffices for our purpose.) We will take as read that the reader is familiar with the standard generic theory of groups of Lie type, as much as is needed in \cite{craven2012un}, so that the unipotent characters of a finite group of Lie type are separated into generic blocks depending only on the particular integer $d$ and not on $\ell$ and $q$. These generic blocks have associated to them a $d$-cuspidal pair $(L,\b\lambda)$, where $L$ is a $d$-split Levi subgroup of $G$ and $\b\lambda$ is a $d$-cuspidal character of $L$.

A set of polynomials in $q$, the \emph{unipotent degrees}, of the non-crystallographic Coxeter groups, were given in \cite{lusztig1993}. These can be separated into generic blocks in the same way as for the groups of Lie type given in \cite{bmm1993}, and this is done in \cite{brouemalle1993}. The labelling of the `unipotent characters' given here is taken from GAP.

We let $G$ be one of the `groups' $I_2(n,q)$, $H_3(q)$ and $H_4(q)$. Write $\Deg(\chi)$ for the generic degree of the unipotent character $\chi$, a polynomial in $q$. Let $B$ be a unipotent block of a $G$ and write $e$ for the number of unipotent characters of $G$ belonging to $B$. Suppose that the `weight' of $B$ is $1$, i.e., the power of $\Phi_d(q)$ dividing $|G|$ is one more than the power dividing each of the generic degrees of the unipotent characters in $B$. The blocks of weight $1$ for groups of Lie type have cyclic defect groups, so have associated Brauer trees. Our goal here is to construct the Brauer trees of the unipotent blocks of weight $1$ for the non-crystallographic groups of Lie type.

\medskip

We now describe the $\pi_{\kappa/d}$-function, where $\kappa$ and $d$ are coprime positive integers (or possibly $d=1$, $\kappa=0$, in which case $\pi_{0/1}(f)=0$ for all $f$). Let $\mc F$ denote the set of all polynomials with real coefficients, whose zeroes are either $0$ or roots of unity. For any $f\in\mc F$, let $a(f)$ be the multiplicity of $0$ as a root of $f$ and $A(f)$ be the degree of $f$; these are often called Lusztig's $a$- and $A$-functions. If $\xi$ is any complex number, write $\Arg_{\kappa/d}(\xi)$ for those $\lambda\in (0,2\pi\kappa/d)$ that are arguments of $\xi$, and if $f\in\mc F$, write $\Arg_{\kappa/d}(f)$ for the multiset consisting of the union of $\Arg_{\kappa/d}(\xi)$ for $\xi$ a non-zero root of $f$ (with multiplicity). Set $\phi_{\kappa/d}(f)$ to be the sum of $|\Arg_{\kappa/d}(f)|$ and half the multiplicity of $1$ as a root of $f$. Finally, set
\[\pi_{\kappa/d}(f)=\frac{\kappa}{d}(a(f)+A(f))+\phi_{\kappa/d}(f).\]
If $\chi$ is a unipotent character in a block with $d$-cuspidal pair $(L,\b\lambda)$, write
\[ \pi_{\kappa/d}(\chi)=\pi_{\kappa/d}(\Deg(\chi))-\pi_{\kappa/d}(\Deg(\b\lambda)).\]

The function $\pi_{\kappa/d}(-)$ conjecturally describes the degree in the cohomology (with non-compact support) over $\bar\Q_\ell$ in which a given unipotent character appears, in a particular  Deligne--Lusztig variety, whose cohomology is meant to yield a derived equivalence between a unipotent block and its Brauer correspondent, consistent with the geometric version of Brou\'e's abelian defect group conjecture. (See \cite{brouemalle1993} for the geometric version of Brou\'e's conjecture, \cite{cravenrouquier2010un} for an introduction to the use of perverse equivalences to this aim, and \cite{craven2012un} for the introduction and analysis of the function $\pi_{\kappa/d}(-)$.)

\medskip

The groups $H_3$ and $H_4$ require $\sqrt 5$ to be present in the field, in which case the cyclotomic polynomials $\Phi_5$, $\Phi_{10}$, $\Phi_{15}$, $\Phi_{20}$ and $\Phi_{30}$ split into two. We always write $\Phi_{5n}''$ for the factor that has $\zeta_{5n}=\e^{2\pi\I/5n}$ as a zero; specifically, we get the following polynomials:
\begin{align*}
\Phi_5'(q)&=q^2+\frac{(1+\sqrt 5)}2 q+1,& \Phi_5''(q)&=q^2+\frac{(1-\sqrt 5)}2 q+1.
\\ \Phi_{10}'(q)&=q^2-\frac{(1-\sqrt 5)}2 q+1,& \Phi_{10}''(q)&=q^2-\frac{(1+\sqrt 5)}2 q+1.
\\ \Phi_{15}'(q)&= q^4+\frac{(1-\sqrt 5)}2(-q^3+q^2-q)+1, &\Phi_{15}''(q)&=q^4+\frac{(1+\sqrt5)}2(-q^3+q^2-q)+1.
\\ \Phi_{20}'(q)&=q^4-\frac{(1-\sqrt5)}2q^2+1&\Phi_{20}''(q)&=q^4-\frac{(1+\sqrt5)}2q^2+1.
\\ \Phi_{30}'(q)&=q^4-\frac{(1-\sqrt5)}2q^2+1,&\Phi_{30}''(q)&=q^4-\frac{(1+\sqrt5)}2q^2+1
\end{align*}
The groups $I_2(n)$ are defined over the field $\Q(\eta+\eta^{-1})$, where $\eta=\zeta_n$. The factors of $x^n-1$ over this field are $\Phi_1(x)$, $\Phi_2(x)$ (if $n$ is even), and
\[ \Phi_n^{(i)}(x)=(x-\eta^i)(x-\eta^{-i}),\]
where $1\leq i<n/2$. We will sometimes use $\Phi$ for some cyclotomic polynomial, over some field of definition.

The order of a finite group of Lie type $G=G(q)$ is given by
\[ q^N\prod_{i=1}^r (q^{a_i}-1),\]
where $N$ is the number of reflections and $a_1,\dots,a_r$ are the degrees of the corresponding Weyl group. These concepts make sense for $I_2$, $H_3$ and $H_4$, and so
\begin{align*} |H_3(q)|&=q^{15}(q^{10}-1)(q^6-1)(q^2-1)=q^{15}\Phi_1^3\Phi_2^3\Phi_3\Phi_5\Phi_6\Phi_{10};
\\ |H_4(q)|&=q^{60}(q^{30}-1)(q^{20}-1)(q^{12}-1)(q^2-1)=q^{60}\Phi_1^4\Phi_2^4\Phi_3^2\Phi_4^2\Phi_5^2\Phi_6^2\Phi_{10}^2\Phi_{12}\Phi_{15}\Phi_{20}\Phi_{30};
\\ |I_2(n,q)|&=q^n(q^n-1)(q^2-1)=q^n\Phi_1^2\Phi_2\prod_{2\leq d\mid n} \Phi_d.
\end{align*}

\section{The combinatorial Brou\'e conjecture}
\label{sec:combbroue}

The full details of the combinatorial Brou\'e conjecture for blocks with cyclic defect groups are given in \cite{craven2012un}, and we provide a brief summary here, sufficient for understanding the rest of the article.

Let $B$ be a unipotent block with $d$-cuspidal pair $(L,\b\lambda)$, and let $B'$ denote its Brauer correspondent, a block of $N_G(D)$, where $D$ is a defect group of $B$. We will construct a bijection between the simple $B$- and $B'$-modules, which will go via the unipotent ordinary $B$-characters and roots of unity in the complex plane.

To $B$ we may associate a \emph{cyclotomic Hecke algebra} $\mc H$, a deformation of the group algebra $\C W$, where $W$ is the \emph{relative Weyl group} associated to $B$. If $B$ has cyclic defect groups, and $e$ simple modules, then $W=Z_e$ is cyclic, and there are $e$ parameters for $\mc H$, each of the form $\omega_i q^{a_i}$, where $q$ is an indeterminate, $\omega_i$ is a root of unity and $a_i$ is a semi-integer.

Since the decomposition matrix of $B$ is (conjecturally) lower unitriangular, there is a one-to-one correspondence between the simple $B$-modules and the unipotent characters belonging to $B$. If $\chi$ is such a unipotent character, with generic degree $f(q)$, write $A(\chi)$ for the degree of $f$, $a(\chi)$ for the multiplicity of $0$ as a root of $f$, and $\aA(\chi)$ for the difference between $a(\chi)+A(\chi)$ and $a(\b\lambda)+A(\b\lambda)$.

To each unipotent character $\chi$ one associates a parameter $\omega_\chi q^{-\aA(\chi)/e}$ of $\mc H$. The root of unity $\omega_\chi$ is related to the \emph{eigenvalue of Frobenius}, which for unipotent characters is $+1$ if the character is $\phi_{a,b}$, and if the character is $H[\alpha]$ for some group of Lie type $H=H(q)$ and root of unity $\alpha$, then $\alpha$ is the eigenvalue of Frobenius. If $d=e$ then $\omega_\chi$ is the eigenvalue of Frobenius, and in other cases it is an appropriate possibly fractional power, depending on the ratio $e/d$. This also applies for the case of $\Phi_d''$, but for $\Phi_d'$ we need to take a power of the root for the $\Phi_d''$-case, see Section \ref{I2n} for example. For a complete list of the $\omega_\chi$, see \cite{brouemalle1993} for classical and most exceptional groups, and the appendix to \cite{craven2012un}, available on the author's webpage, for all unipotent blocks with cyclic defect groups of exceptional groups. If $u_1,\dots,u_e$ are the parameters of $\mc H$, then the \emph{relative degree} associated to $u_i$ is
\[ \prod_{j\neq i} \frac{u_j}{u_i-u_j}.\]
Up to a global scaling function, the relative degree associated to $u_i$ is simply the generic degree of the unipotent character associated to $u_i$.

For the combinatorial Brou\'e conjecture, we need to renormalize the parameters of the cyclotomic Hecke algebra, by taking inverses and complex conjugates, hence associating to each $\chi$ the quantity $\omega_\chi q^{\aA(\chi)/e}$. The map $q\mapsto \e^{2\pi\I\kappa/d}$ induces a bijection between the ordinary characters $\chi$ and various roots of unity (normally all $e$th roots, sometimes all $2e$th roots that are not $e$th roots) via the quantity $\omega_\chi q^{\aA(\chi)/e}$. Hence the composition of these two bijections associates, to each simple $B$-module, a specific root of unity.

Choose a unipotent character $\psi$ that has minimal $\pi_{\kappa/d}$-value among the unipotent characters in the block $B$ (there might be more than one); then $\psi$ reduces modulo $\ell$ to a simple $B$-module, whose Green correspondent $M$ is either simple or is the Heller translate of a simple module. (This is not true for arbitrary blocks with cyclic defect group, and in general it depends on a certain endopermutation module. For groups of Lie type, $\Norm_G(D)=\Norm_G(\Omega_1(D))$, where $D$ is the defect group and $\Omega_1(D)$ is the subgroup of order $p$ of $D$, in which case the statement does hold. See \cite{linckelmann1996} for full details.)

The Brauer tree of $B'$, the Brauer correspondent of $B$, is a star with $e$ edges; embed this star into the complex plane so that the exceptional node is at $0$, and the edges are equally spaced out around the exceptional node with the vertices lying at complex numbers of modulus $1$. Specifying the position of one of the edges will associate to every simple $B'$-module a particular complex number of modulus $1$. Place $M$ at the position $\omega_\chi$ in the complex plane: if $M$ is simple then this fixes the rotation of the Brauer tree of $B'$ and so associates roots of unity to each simple $B'$-module; otherwise the characters corresponding to $\Omega(M)$ and $\Omega^{-1}(M)$ should be placed so that $\omega_\chi$ is half-way between them, and again this associates to each simple $B'$-module a particular root of unity. Hence we have constructed a bijection between the simple $B'$-modules and various roots of unity, and composing with the previous bijection produces a bijection between the simple $B$-modules and the simple modules of its Brauer correspondent $B'$.

In \cite{chuangrouquierun} the theory of perverse equivalences is developed (see also \cite{rouquier2006}, \cite{cravenrouquier2010un}, and particularly \cite{craven2012un}):  a perverse equivalence is a special type of derived equivalence between $B$ and $B'$, and can be combinatorially described via a stable equivalence between $B$ and $B'$, a perversity function $\pi(-)$ on the simple $B$-modules, and a bijection between the simple $B$- and $B'$-modules. In \cite{craven2012un} it is proved that, with the possible exception of two unipotent blocks for $E_8$ and one for ${}^2\!F_4$, whenever the block $B$ has cyclic defect groups there is a perverse equivalence between $B$ and $B'$ with perversity function $\pi_{\kappa/d}(-)$, the bijection given above, and with induction and restriction as the stable equivalence. This statement is the \emph{combinatorial Brou\'e conjecture}. One particular aspect of the perversity function is that it increases towards the exceptional node, which is used in Section \ref{I2n}.

Furthermore, given a function $f(-)$ on the simple $B'$-modules, there is a unique block (up to Morita equivalence) $C$ with a perverse equivalence between $C$ and $B'$ with some bijection between the simple $C$- and $B'$-modules such that the perversity function, pulled through the bijection, yields $f$. In other words, since the bijection between the simple modules and perversity function are both determined by the cyclotomic Hecke algebra associated to the block, assuming the combinatorial Brou\'e conjecture for the block $B$ one may reconstruct the Brauer tree of $B$.

All unipotent blocks of $I_2(n,q)$, $H_3(q)$ and $H_4(q)$ have associated cyclotomic Hecke algebras (see \cite[8.3]{brouemalle1993}, although notably some cases are missing there, which are dealt with here), and so we may assume the combinatorial Brou\'e conjecture for these `groups' and reconstruct the Brauer tree. We can also use standard arguments given in the next section to reconstruct, at least partially, the Brauer tree. In all cases the tree constructed from the combinatorial Brou\'e conjecture is entirely consistent with the information arising out of the standard arguments, implying that the representation theory of these `groups' still makes sense, even though the objects themselves do not.

\section{Arguments for Brauer trees}
\label{sec:Brauertreeargs}
In the past, various arguments have been used to understand the Brauer trees of groups of Lie type, and many of these can be carried over to the non-crystallographic case. We summarize these arguments now:
\begin{enumerate}
\item Parity argument. In a block with cyclic defect group, the sum of two characters that are labelled by adjacent vertices on the Brauer tree is a projective character, so has degree divisible by $\ell$; this therefore partitions the set of characters into two. In particular, for groups of Lie type the generic degree $\Deg(\chi)/\Deg(\b\lambda)$ is congruent to either $\alpha$ or $-\alpha$ modulo $\Phi$, where $\alpha$ is some positive real number. We refer to the two sets of characters as $+$-type and $-$-type; the argument is that two $+$-type characters cannot be adjacent, and similarly two $-$-type characters cannot be adjacent.

\item Degree argument. The degree of an irreducible non-exceptional character $\chi$ is the sum of the dimensions of the simple modules labelled by edges incident to $\chi$ on the Brauer tree; since dimensions of simple modules must be positive, this places constraints on the possible configurations. Broadly speaking, the degree of a character, as a polynomial in $q$, must increase towards the exceptional node. Sometimes (for $H_3$ and $d=5'$, $6$, $10''$ and $H_4$, $d=10'$) we require $q>2$ for the degree argument to work, simply because $q=2$ is too small for one polynomial to be larger than another.

\item M\"uller argument. The Brauer trees of the Hecke algebras of Coxeter groups have been determined by Geck \cite{geck1992} for Weyl groups, and by M\"uller \cite{muller1997} for the non-crystallographic case, and are always lines. This line consists of the principal-series characters, and forms a subtree of the Brauer tree of the block, connected to the exceptional node.

\item The real stem. Suppose that $B$ is a real block. The subset of real characters forms a subtree of the Brauer tree that is a line with the exceptional node somewhere along its length. By the M\"uller argument above, we already know one side of the exceptional node, and all non-principal-series real characters form the other side. Together with the degree argument, this is enough to determine the complete real stem. Hence it is only the non-real characters in real blocks that need to be placed.

\item Harish-Chandra induction argument. We can use Harish-Chandra induction, which is defined for the unipotent characters of $I_2(n)$, $H_3$ and $H_4$, to induce projective characters for a subgroup, and get projective characters for the overgroup. This will not work for cuspidal unipotent characters, but often allows us to understand the location of the $I_2(5)$-series characters for $H_3(q)$ and $H_4(q)$. Harish-Chandra induction sends projectives to projectives, and indeed the same is true of Harish-Chandra restriction, which will be needed for $H_3$.)
\end{enumerate}

These arguments can produce many decomposition numbers, but not the planar embedding. Alvis--Curtis duality, which can be used to get more information about the Brauer trees for genuine groups of Lie type, is not used in this paper, but see \cite{cdr2012un} for an indication of how it can be used.

Finally, when $\Phi$ is the `Coxeter polynomial', i.e., the order of the Coxeter torus, so $\Phi_n^{(1)}$ for $I_2(n,q)$, $\Phi_{10}''$ for $H_3(q)$ and $\Phi_{30}''$ for $H_4(q)$, in \cite{hlm1995} Hiss, L\"ubeck and Malle provide a conjecture for the planar-embedded Brauer tree for genuine groups of Lie type, finally proved in the remaining cases by Dudas and Rouquier \cite{dudasrouquier2012un}. This states that the tree consists of lines emanating from the exceptional node, with each line containing characters with the same eigenvalue of Frobenius (a root of unity), and ordered around the exceptional node in the same manner as the eigenvalues of Frobenius are arranged around $0$ on the complex plane, in increasing argument (as a complex number). As part of our work here, we will see that the Hiss--L\"ubeck--Malle conjecture is compatible with the combinatorial Brou\'e conjecture, so that it extends to the non-crystallographic case.

In \cite{cdr2012un} the Deligne--Lusztig variety associated to the Coxeter torus is studied at primes \emph{not} dividing the Coxeter polynomial: since the variety does not exist for $H_3$, $H_4$ and $I_2(n)$, the method cannot be used, but the statement can still be translated over, and this could also be used to determine more Brauer trees. However, we do not pursue this idea here.

\section{The group $I_2(n,q)$}
\label{I2n}
Let $\eta=\e^{2\pi\I/n}$. Since the $p'$-part of the order of $I_2(n,q)$ is $\Phi_1\Phi_2(q^n-1)$, the possible irreducible factors (over the field $\Q(\eta+\eta^{-1})$) that $\ell$ can divide are $\Phi_1(q)$, $\Phi_2(q)$ and $\Phi_n^{(i)}(q)$ for $1\leq i<n/2$.

Depending on whether $n$ is odd or even, we get different unipotent degrees. In the following table, the last two rows only apply if $n$ is even. We include the eigenvalue of Frobenius.
\begin{center}\begin{tabular}{lcc}
\hline Character & Degree & Eigenvalue
\\\hline $\phi_{1,0}$ & $1$ & $1$
\\ $\phi_{1,n}$ & $q^n$ & $1$
\\ $\phi_{2,i}$, $i<n/2$ & $\D\frac{(1-\eta^i)(1-\eta^{-i})}{n}\frac{q\Phi_2(q^n-1)}{\Phi_1 \Phi_n^{(i)}}$ & $1$
\\ $I_2(n)[i,j]$, $i<j<n-i$ & $\D\frac{\eta^i+\eta^{-i}-\eta^j-\eta^{-j}}n \frac{q(q^2-1)(q^n-1)}{\Phi_n^{(i)}\Phi_n^{(j)}}$ & $\eta^{-ij}$
\\\ $\phi_{1,n/2}'$ & $\D\frac{2}{n}\frac{q(q^n-1)}{\Phi_1\Phi_2}$ & $1$
\\ $\phi_{1,n/2}''$ & $\D\frac{2}{n}\frac{q(q^n-1)}{\Phi_1\Phi_2}$ & $1$
\\\hline\end{tabular}\end{center}
In Lusztig's notation from \cite{lusztig1993}, the character $\phi_{1,0}$ is called $1$, the character $\phi_{1,n}$ is called $S$, the characters $\phi_{2,i}$ are $\rho_i$, the characters $I_2(n)[i,j]$ for $j<n/2$ are called $\rho_{i,j}$, and for $j>n/2$ are called $\rho_{i,n-j}'$. If $n$ is even, then $\phi_{1,n/2}'$ and $\phi_{1,n/2}''$ are $\ep'$ and $\ep''$, and $I_2(5)[i,n/2]$ is $\rho_i'$.

We make a few remarks about the characters $I_2(n)[i,j]$. Firstly, the dual of $I_2(n)[i,j]$ is $I_2(n)[i,n-j]$. Secondly, it is often useful to allow ourselves to swap the indices $i$ and $j$ (so that $I_2(n)[j,i]=I_2(n)[i,j]$), and also simultaneously replace $i$ and $j$ by $n-i$ and $n-j$ (so that, if $n=10$ for example, the pairs $[1,4]$, $[4,1]$, $[9,6]$ and $[6,9]$ all label the same character). We do this purely to make writing down formulae much easier. Notice that even with these extra ways of defining a character, each character is still well defined.

There is no unipotent character with a single power of $\Phi_1$ dividing it, and so there are no unipotent $\Phi_1$-blocks with cyclic defect group. The exact power of all other cyclotomic factors is $1$ (unless $d=2$ and $n$ is even, in which case there are no unipotent blocks of weight $1$), and so the principal block is of weight $1$ for all other factors, and indeed is the only unipotent block.

Let $1\leq i<n/2$, write $d'=(n,i)$, and consider the principal $\Phi_n^{(i)}$-block, so that $d=n/d'$ and $\kappa=i/d',(n-i)/d'$ (plus any multiple of $d$). We claim that the parameters of the cyclotomic Hecke algebra of the (only) unipotent block are $1,q,q^2$, and $\eta^jq$ for $i\neq j\leq n/2$, with the table below giving the $\pi_{\kappa/d}$-function for $\kappa=i/d'$ and $\kappa=(n-i)/d'$.
\begin{center}\begin{tabular}{lcccc}
\hline Character & $A(-)$ & Parameter & $\kappa=i/d'$ & $\kappa=(n-i)/d'$
\\\hline $\phi_{1,0}$ & $0$ & $q^0$ & $0$ & $0$
\\ $\phi_{2,i}$ & $n-1$ & $q$ & $2i-1$ & $2(n-i)-1$
\\ $\phi_{1,n}$ & $n$ & $q^2$ & $2i$ & $2(n-i)$
\\ $I_2(n)[i,j]$ & $n-1$ & $\eta^{-j}q$ & $2i-\begin{cases}1 & j<i\\ 0& j>i\end{cases}$ & $2(n-i)-\begin{cases}1 & j<i\\ 2& j>i\end{cases}$
\\ $I_2(n)[i,n-j]$ & $n-1$ & $\eta^jq$ & $2i-\begin{cases}1 & j<i\\ 0& j>i\end{cases}$ & $2(n-i)-\begin{cases}1 & j<i\\ 2& j>i\end{cases}$
\\\hline\end{tabular}\end{center}
The real stem is easy to construct, since by \cite{muller1997}, the principal-series characters are $\phi_{1,0}$, connected to $\phi_{2,i}$, connected to $\phi_{1,n}$, with the exceptional node connected to that. If $n$ is even then there is one more real character, namely $I_2(n)[i,n/2]$, which must be located at the other side of the exceptional. The characters $I_2(n)[i,j]$ and $I_2(n)[i,n-j]$ with $i<j<n/2$ must be connected to $\phi_{1,n}$ by a parity argument, and this leaves $I_2(n)[i,j]$ and $I_2(n)[i,n-j]$ with $j<i<n/2$, which by parity must be connected to the exceptional or $\phi_{2,i}$. For some values of $i$ and $j$, twice the degree of $I_2(n)[i,j]$ is less than that of $\phi_{2,i}$, and sometimes it is greater. Thus it is difficult to know where to place these non-real characters on the Brauer tree. If we assume the combinatorial Brou\'e conjecture, however, then it is clear where they should go, since $\pi_{\kappa/d}(I_2(n)[i,j])=\pi_{\kappa/d}(\phi_{2,i})$ for the relevant $j$. (Indeed, this proves that we do not need the full combinatorial Brou\'e conjecture here, but merely the fact that the $\pi_{\kappa/d}$-function increases towards the exceptional node, as mentioned in Section \ref{sec:combbroue}.)

The planar embedding, which cannot be recovered without assuming the combinatorial Brou\'e conjecture, can be easily understood using the cyclotomic Hecke algebra of the block, whose parameters we gave above, and is the following.
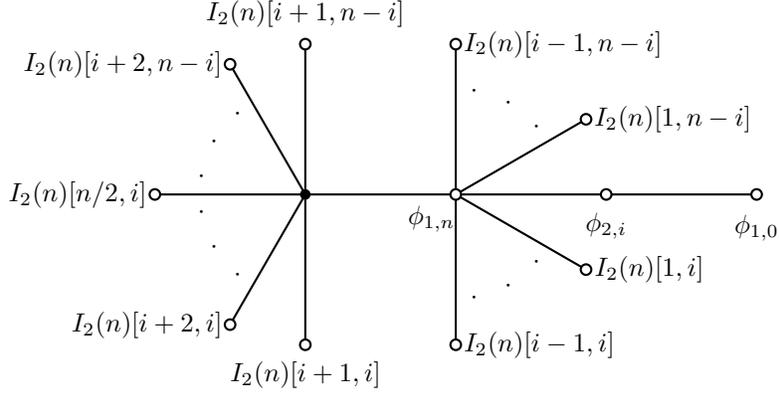
\begin{center}\begin{tikzpicture}[thick,scale=2]

\draw (3,1.2) node {${I_2(n)[i+1,n-i]}$};
\draw (3,-1.2) node {${I_2(n)[i+1,i]}$};

\draw (2,0) -- (6,0);
\draw (3,1) -- (3,-1);
\draw (3,0) -- (2.5,0.866);
\draw (3,0) -- (2.5,-0.866);
\draw (4,1) -- (4,-1);
\draw (4,0) -- (4.866,0.5);
\draw (4,0) -- (4.866,-0.5);

\begin{scope}[shift={(4,0)}]
\foreach \i in {1,...,3} \node (\i) at (20+20*\i:0.7) {$\cdot$};
\foreach \i in {1,...,3} \node (\i) at (-20-20*\i:0.7) {$\cdot$};
\end{scope}
\begin{scope}[shift={(3,0)}]
\foreach \i in {1,...,3} \node (\i) at (110+20*\i:0.7) {$\cdot$};
\foreach \i in {1,...,3} \node (\i) at (-110-20*\i:0.7) {$\cdot$};
\end{scope}
\draw (3,0) node [fill=black!100] (ld) {};
\draw (4,0) node [draw,label=below left:$\phi_{1,n}$] (l4) {};
\draw (5,0) node [draw,label=below:$\phi_{2,i}$] (l4) {};
\draw (6,0) node [draw,label=below:$\phi_{1,0}$] (l4) {};

\draw (4.866,0.5) node [draw,label=right:${I_2(n)[1,n-i]}$] (l4) {};
\draw (4.866,-0.5) node [draw,label=right:${I_2(n)[1,i]}$] (l4) {};

\draw (4,1) node [draw,label=right:${I_2(n)[i-1,n-i]}$] (l4) {};
\draw (4,-1) node [draw,label=right:${I_2(n)[i-1,i]}$] (l4) {};

\draw (2.5,0.866) node [draw,label=left:${I_2(n)[i+2,n-i]}$] (l4) {};
\draw (2.5,-0.866) node [draw,label=left:${I_2(n)[i+2,i]}$] (l4) {};

\draw (3,1) node [draw] (l4) {};
\draw (3,-1) node [draw] (l4) {};

\draw (2,0) node [draw,label=left:${I_2(n)[n/2,i]}$] (l4) {};

\end{tikzpicture}\end{center}
(The node $I_2(n)[n/2,i]$ only exists if $n$ is even, and lies on the real stem in this case.) This of course encapsulates the cases of all Brauer trees for $\GL_3(q)$, $\Sp_4(q)$ and $G_2(q)$.

If $i=1$, so that $\Phi_n^{(1)}$ is the Coxeter torus, then all of the cuspidal unipotent characters are arranged around the exceptional node in order of increasing argument of eigenvalue of Frobenius, and we see that the Hiss--L\"ubeck--Malle conjecture is verified (see Section \ref{sec:Brauertreeargs}).

Finally, if $n$ is odd then there is a principal block with cyclic defect groups for $d=2$, consisting of $\phi_{1,0}$ and $\phi_{1,n}$, so the Brauer tree is a line with the exceptional connected to $\phi_{1,n}$.

\section{The group $H_3(q)$}

For this group there are two blocks of weight $1$ for $d=1$, each with two unipotent characters, the same for $d=2$, and the principal block for $d=3,5',5'',6,10',10''$. We first deal with the principal blocks, and then with $d\leq 2$. In the case where $d=3$, all unipotent characters lie in the principal series, and so from \cite{muller1997} we already know almost all of the decomposition matrix. Using the standard degree argument we can tell which of the ordinary characters of the two blocks of the Hecke algebra are connected to the exceptional node, and end with the following tree.
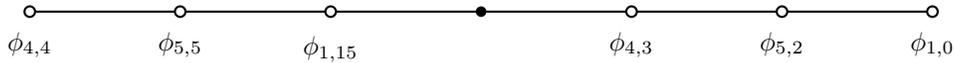
\begin{center}\begin{tikzpicture}[thick,scale=2]
\draw (0,0) -- (6,0);
\draw (0,0) node [draw,label=below:$\phi_{4,4}$] (l0) {};
\draw (1,0) node [draw,label=below:$\phi_{5,5}$] (l1) {};
\draw (2,0) node [draw,label=below:$\phi_{1,15}$] (l2) {};
\draw (3,0) node [fill=black!100] (ld) {};
\draw (4,0) node [draw,label=below:$\phi_{4,3}$] (l4) {};
\draw (5,0) node [draw,label=below:$\phi_{5,2}$] (l4) {};
\draw (6,0) node [draw,label=below:$\phi_{1,0}$] (l4) {};
\end{tikzpicture}\end{center}

For $\Phi_5''$, we can use the tables in \cite{muller1997}, i.e., a M\"uller argument, together with a degree argument, to get the real stem, consisting of six characters. The remaining four non-real characters -- $I_2(5)[1,3];1$, $I_2(5)[1,2];1$, $I_2(5)[1,3];\ep$ and $I_2(5)[1,2];\ep$ -- are of $+$-type, so either connected to the exceptional node or to $\phi_{4,3}$ or $\phi_{4,4}$; a degree argument proves that $I_2(5)[1,2];\ep$ and its conjugate must be connected to the exceptional node. In order to prove that the other $I_2(5)$-characters are connected to the exceptional node we need to determine the tree for $H_4$ and $d=5''$, which we do in the next section. To get the full planar embedding however, we have to assume the combinatorial Brou\'e conjecture; doing this yields the following tree.
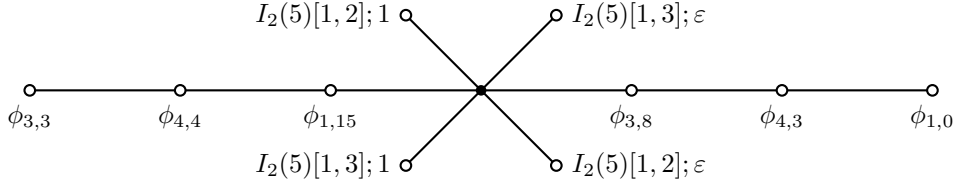
\begin{center}\begin{tikzpicture}[thick,scale=2]
\tikzstyle{every node}=[rectangle]

\draw (4.05,0.5) node {$I_2(5)[1,2];\ep$};
\draw (1.95,1.5) node {$I_2(5)[1,2];1$};
\draw (4.05,1.5) node {$I_2(5)[1,3];\ep$};
\draw (1.95,0.5) node {$I_2(5)[1,3];1$};

\draw (6,0.82) node {$\phi_{1,0}$};
\draw (5,0.82) node {$\phi_{4,3}$};
\draw (4,0.82) node {$\phi_{3,8}$};
\draw (2,0.82) node {$\phi_{1,15}$};
\draw (1,0.82) node {$\phi_{4,4}$};
\draw (0,0.82) node {$\phi_{3,3}$};

\tikzstyle{every node}=[circle, fill=black!0, inner sep=0pt, minimum width=4pt]

\draw (0,1) -- (6,1);
\draw (3.5,0.5) -- (2.5,1.5);
\draw (2.5,0.5) -- (3.5,1.5);

\draw (6,1) node [draw] (l0) {};
\draw (5,1) node [draw] (l1) {};
\draw (4,1) node [draw] (l1) {};
\draw (0,1) node [draw] (l1) {};
\draw (3,1) node [fill=black!100] (ld) {};
\draw (2,1) node [draw] (l0) {};
\draw (1,1) node [draw] (l1) {};
\draw (2.5,0.5) node [draw] (l0) {};
\draw (2.5,1.5) node [draw] (l1) {};
\draw (3.5,0.5) node [draw] (l0) {};
\draw (3.5,1.5) node [draw] (l1) {};
\end{tikzpicture}\end{center}

When $d=5'$, there are again ten unipotent characters in the principal block, six of which come from the principal series; as $5'$ is the algebraic conjugate of $5''$, we must use the algebraic conjugates of the tables in \cite{muller1997}, which means that $\phi_{1,0}$, $\phi_{4,3}$, and the algebraic conjugate of $\phi_{3,8}$, namely $\phi_{3,6}$, lie in one block, and $\phi_{3,1}$, $\phi_{4,4}$ and $\phi_{1,15}$ lie in the other; a degree argument fixes the shape of the real stem. The remaining four characters are all of $-$-type, so cannot be connected to the exceptional node. Moreover, $I_2(5)[1,2];\ep$ has larger degree than $\phi_{3,6}$ for $q>2$ (the degree of the polynomial is the same, but the leading coefficient is larger) and so $I_2(5)[1,2];\ep$ and its dual must be connected to $\phi_{1,15}$ by a degree argument. As with $d=5''$, in order to prove that the other $I_2(5)$-characters are connected to $\phi_{3,6}$ we need to determine the tree for $H_4$ and $d=5'$, which we do in the next section.

Standard arguments cannot determine the planar embedding, and the combinatorial Brou\'e conjecture predicts that it looks like the following.
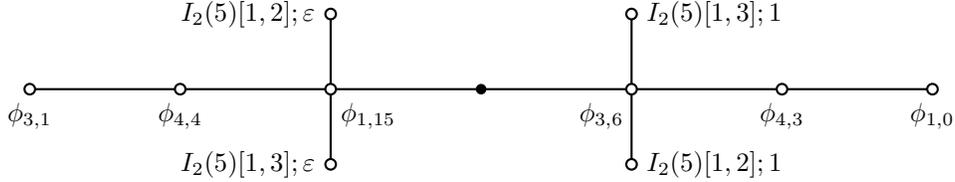
\begin{center}\begin{tikzpicture}[thick,scale=2]
\tikzstyle{every node}=[rectangle]

\draw (4.55,0.5) node {$I_2(5)[1,2];1$};
\draw (1.45,1.5) node {$I_2(5)[1,2];\ep$};
\draw (4.55,1.5) node {$I_2(5)[1,3];1$};
\draw (1.45,0.5) node {$I_2(5)[1,3];\ep$};

\draw (6,0.82) node {$\phi_{1,0}$};
\draw (5,0.82) node {$\phi_{4,3}$};
\draw (3.8,0.82) node {$\phi_{3,6}$};
\draw (2.25,0.82) node {$\phi_{1,15}$};
\draw (1,0.82) node {$\phi_{4,4}$};
\draw (0,0.82) node {$\phi_{3,1}$};

\tikzstyle{every node}=[circle, fill=black!0, inner sep=0pt, minimum width=4pt]

\draw (0,1) -- (6,1);
\draw (4,0.5) -- (4,1.5);
\draw (2,0.5) -- (2,1.5);

\draw (6,1) node [draw] (l0) {};
\draw (5,1) node [draw] (l1) {};
\draw (4,1) node [draw] (l1) {};
\draw (0,1) node [draw] (l1) {};
\draw (3,1) node [fill=black!100] (ld) {};
\draw (2,1) node [draw] (l0) {};
\draw (1,1) node [draw] (l1) {};
\draw (2,0.5) node [draw] (l0) {};
\draw (2,1.5) node [draw] (l1) {};
\draw (4,0.5) node [draw] (l0) {};
\draw (4,1.5) node [draw] (l1) {};
\end{tikzpicture}\end{center}
The relationship between $H_2=I_2(5)$ and $H_3$ is very similar to the relationship between $E_6$ and $E_7$; when $d$ is odd (including $d=5'$ and $d=5''$), the Brauer tree for $H_3$ is a doubling of that of $I_2(5)$, just as the Brauer tree for $E_7$ is a doubling of that of $E_6$, for $d$ odd ($d=3,5,9$). As of yet though, no method of proof exists that reconstructs this doubling directly for $E_7$; should this be achieved, it might be possible to determine the planar embedding for $d=5''$ in the same way.

\medskip

We now move on to $d=6$. Here the real stem is given by \cite{muller1997}, so it remains to locate the characters $H_3[\pm \I]$. However, they are of $+$-type, so can either be connected to $\phi_{1,15}$ or $\phi_{5,2}$, the two $-$-type characters. A degree argument shows that they cannot be connected to $\phi_{5,2}$ for $q>2$, so we get the following diagram.
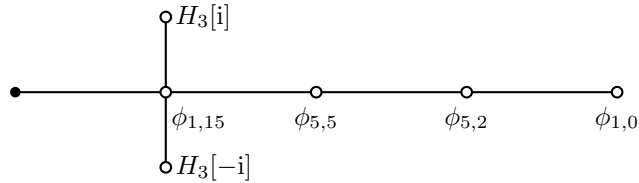
\begin{center}\begin{tikzpicture}[thick,scale=2]
\tikzstyle{every node}=[rectangle]

\draw (1.25,1.5) node {$H_3[\I]$};
\draw (1.32,0.5) node {$H_3[-\I]$};

\draw (4,0.82) node {$\phi_{1,0}$};
\draw (3,0.82) node {$\phi_{5,2}$};
\draw (2,0.82) node {$\phi_{5,5}$};
\draw (1.22,0.82) node {$\phi_{1,15}$};
\tikzstyle{every node}=[circle, fill=black!0, inner sep=0pt, minimum width=4pt]

\draw (0,1) -- (4,1);
\draw (1,0.5) -- (1,1.5);

\draw (4,1) node [draw] (l0) {};
\draw (3,1) node [draw] (l1) {};
\draw (2,1) node [draw] (l1) {};
\draw (1,1) node [draw] (l1) {};
\draw (0,1) node [fill=black!100] (ld) {};
\draw (1,0.5) node [draw] (l0) {};
\draw (1,1.5) node [draw] (l1) {};
\end{tikzpicture}\end{center}

The penultimate case is that of $d=10''$, the Coxeter case. Here the Hiss--L\"ubeck--Malle conjecture, detailed in Section \ref{sec:Brauertreeargs}, suggests the planar-embedded Brauer tree. We get the real stem from \cite{muller1997}, and know that $I_2(5)[1,3];1+I_2(5)[1,3];\ep$ is a projective character by Harish-Chandra induction of $I_2(5)[1,3]$ (which is projective) from the $I_2(5)$ subgroup, so that $I_2(5)[1,3];1$ and $I_2(5)[1,3];\ep$ are connected. By degree and parity arguments, these $I_2(5)$-series characters must be connected to the exceptional node, only for $q>2$. However, using only degree and parity arguments, $H_3[\pm\I]$ could be connected to one of two characters: the exceptional or $\phi_{3,6}$. The one that follows the guideline of Hiss, L\"ubeck and Malle is that $H_3[\pm\I]$ are connected to the exceptional node, and this is also the result of applying the combinatorial Brou\'e conjecture.

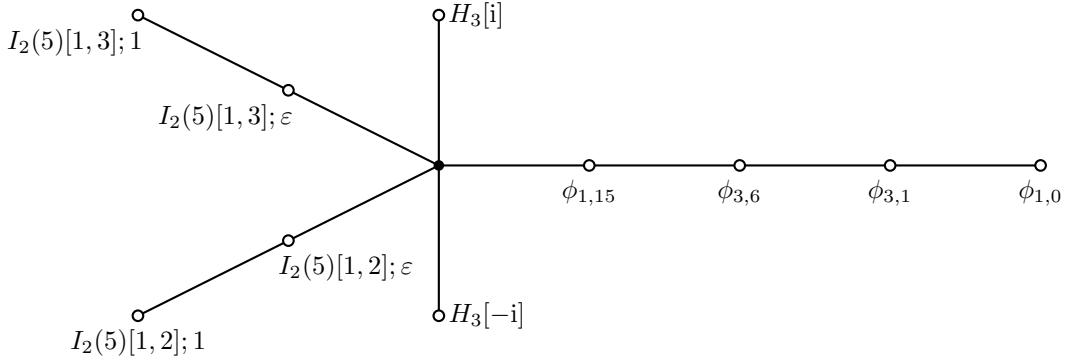
\begin{center}\begin{tikzpicture}[thick,scale=2]
\tikzstyle{every node}=[rectangle]

\draw (4,0.82) node {$\phi_{1,0}$};
\draw (3,0.82) node {$\phi_{3,1}$};
\draw (2,0.82) node {$\phi_{3,6}$};
\draw (1,0.82) node {$\phi_{1,15}$};

\draw (0.25,2) node {$H_3[\I]$};
\draw (0.32,0) node {$H_3[-\I]$};

\draw (-2.42,1.82) node {$I_2(5)[1,3];1$};
\draw (-1.42,1.32) node {$I_2(5)[1,3];\ep$};

\draw (-2,-0.18) node {$I_2(5)[1,2];1$};
\draw (-0.62,0.32) node {$I_2(5)[1,2];\ep$};

\tikzstyle{every node}=[circle, fill=black!0, inner sep=0pt, minimum width=4pt]

\draw (0,1) -- (4,1);
\draw (0,0) -- (0,2);
\draw (0,1) -- (-2,2);
\draw (0,1) -- (-2,0);

\draw (4,1) node [draw] (l0) {};
\draw (3,1) node [draw] (l1) {};
\draw (2,1) node [draw] (l1) {};
\draw (1,1) node [draw] (l1) {};
\draw (0,1) node [fill=black!100] (ld) {};
\draw (0,0) node [draw] (l0) {};
\draw (0,2) node [draw] (l1) {};
\draw (-1,0.5) node [draw] (l0) {};
\draw (-1,1.5) node [draw] (l1) {};
\draw (-2,0) node [draw] (l0) {};
\draw (-2,2) node [draw] (l1) {};
\end{tikzpicture}\end{center}

We are left with the case $d=10'$. In this case again, the structure of the Brauer tree cannot be determined by standard arguments, as the cuspidal characters $H_3[\pm\I]$ could again be in more than one place. The structure of the real stem can be obtained from \cite{muller1997} (taking algebraic conjugates as with $d=5'$), and $I_2(5)[1,3];1$ and $I_2(5)[1,3];\ep$ are again connected by a Harish-Chandra induction argument, with $I_2(5)[1,3];\ep$ being connected to $\phi_{1,15}$ by parity and degree arguments. However, the exact position of $H_3[\pm \I]$ cannot be seen in this way: they could be connected to $I_2(5)[1,2];\ep$ and its conjugate, or to $\phi_{3,8}$.

Thus the following Brauer tree is therefore the result of applying the combinatorial Brou\'e conjecture.
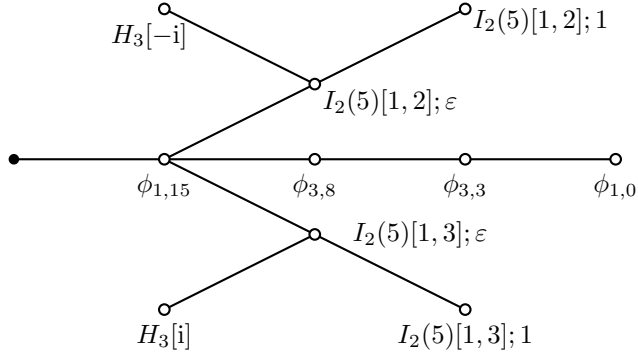
\begin{center}\begin{tikzpicture}[thick,scale=2]

\tikzstyle{every node}=[rectangle]
\tikzstyle{every node}=[circle, fill=black!0, inner sep=0pt, minimum width=4pt]

\draw (0.9,1.82) node {$H_3[-\I]$};
\draw (3.5,1.92) node {$I_2(5)[1,2];1$};
\draw (2.5,1.38) node {$I_2(5)[1,2];\ep$};

\draw (1,-0.18) node {$H_3[\I]$};
\draw (3,-0.18) node {$I_2(5)[1,3];1$};
\draw (2.7,0.5) node {$I_2(5)[1,3];\ep$};

\draw (4,0.82) node {$\phi_{1,0}$};
\draw (3,0.82) node {$\phi_{3,3}$};
\draw (2,0.82) node {$\phi_{3,8}$};
\draw (1,0.82) node {$\phi_{1,15}$};

\tikzstyle{every node}=[circle, fill=black!0, inner sep=0pt, minimum width=4pt]

\draw (0,1) -- (4,1);
\draw (1,1) -- (3,2);
\draw (2,1.5) -- (1,2);
\draw (1,1) -- (3,0);
\draw (2,0.5) -- (1,0);

\draw (4,1) node [draw] (l0) {};
\draw (3,1) node [draw] (l1) {};
\draw (2,1) node [draw] (l1) {};
\draw (1,1) node [draw] (l1) {};
\draw (0,1) node [fill=black!100] (ld) {};

\draw (3,2) node [draw] (l0) {};
\draw (1,2) node [draw] (l1) {};
\draw (2,1.5) node [draw] (l0) {};

\draw (3,0) node [draw] (l1) {};
\draw (1,0) node [draw] (l0) {};
\draw (2,0.5) node [draw] (l1) {};
\end{tikzpicture}\end{center}

We briefly mention the non-principal blocks when $d=1$ and $d=2$. In both cases there are two non-principal unipotent blocks: for $d=1$, the $1$-cuspidal pairs involved are $(\Phi_1.I_2(5),I_2(5)[1,2])$ and its dual, and for $d=2$ the $2$-cuspidal pairs involved are $(\Phi_2.{}^2\!I_2(5),{}^2\!I_2(5)[1,2])$ and its dual. This is needed when calculating the $\pi_{\kappa/d}$-function to test Brou\'e's conjecture but not for computing the Brauer trees: since the two unipotent characters in this block are $I_2(5)[1,2];1$ and $I_2(5)[1,2];\ep$ in both cases, we see by a parity argument that the Brauer trees are lines, with the exceptional in the middle for $d=1$ and at the end, connected to $I_2(5)[1,2];\ep$ for $d=2$.

\section{$H_4(q)$}

As with $H_3$, we first deal with the principal blocks and then the non-principal blocks. There are principal blocks in the cases $d=12$, $15'$, $15''$, $20'$, $20''$, $30'$ and $30''$, and non-principal blocks in the case $d=1$, $2$, $5'$, $5''$, $10'$ and $10''$.

We start with $d=12$. Here, we use \cite{muller1997} and a degree argument to understand the real stem (the character $H_4^\mathrm{VI}[-1]$ is also real and connected to the exceptional node). All non-real characters are of $-$-type, so not connected to the exceptional node, and can be connected to any of $\phi_{1,60}$ (of degree $q^{60}$), or $\phi_{48,9}$ or $H_4^\mathrm{VI}[-1]$ (both with leading term in their degrees $q^{54}/12$). Since all the non-real unipotent characters here have leading terms $q^{54}/4$ (for $H_4[\pm\I]$) and $q^{54}/6$ (for the others), the non-real characters must be attached to $\phi_{1,60}$.

The planar embedding cannot be recovered, as with previous cases, and so we have to consider the cyclotomic Hecke algebra for this. Assuming the combinatorial version of Brou\'e's conjecture holds, we get the following Brauer tree.
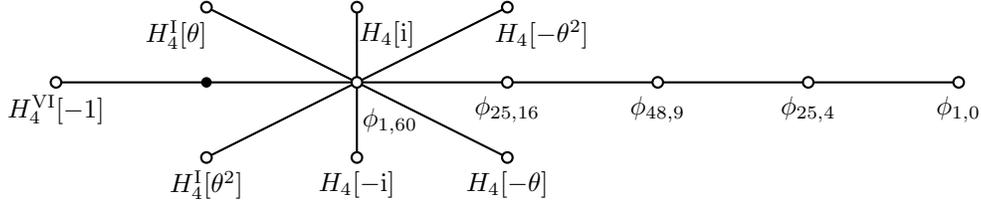
\begin{center}\begin{tikzpicture}[thick,scale=2]
\tikzstyle{every node}=[rectangle]

\draw (5,0.82) node {$\phi_{1,0}$};
\draw (4,0.82) node {$\phi_{25,4}$};
\draw (3,0.82) node {$\phi_{48,9}$};
\draw (2,0.82) node {$\phi_{25,16}$};
\draw (1.22,0.74) node {$\phi_{1,60}$};
\draw (-1,0.82) node {$H_4^\mathrm{VI}[-1]$};

\draw (0,0.32) node {$H_4^\mathrm{I}[\theta^2]$};
\draw (1,0.32) node {$H_4[-\I]$};
\draw (2,0.32) node {$H_4[-\theta]$};

\draw (-0.2,1.32) node {$H_4^\mathrm{I}[\theta]$};
\draw (1.2,1.32) node {$H_4[\I]$};
\draw (2.23,1.32) node {$H_4[-\theta^2]$};

\tikzstyle{every node}=[circle, fill=black!0, inner sep=0pt, minimum width=4pt]

\draw (-1,1) -- (5,1);
\draw (0,0.5) -- (2,1.5);
\draw (2,0.5) -- (0,1.5);
\draw (1,0.5) -- (1,1.5);

\draw (5,1) node [draw] (l0) {};
\draw (4,1) node [draw] (l0) {};
\draw (3,1) node [draw] (l1) {};
\draw (2,1) node [draw] (l1) {};
\draw (1,1) node [draw] (l1) {};
\draw (-1,1) node [draw] (l1) {};
\draw (0,1) node [fill=black!100] (ld) {};
\draw (0,0.5) node [draw] (l0) {};
\draw (1,0.5) node [draw] (l0) {};
\draw (2,0.5) node [draw] (l1) {};
\draw (0,1.5) node [draw] (l0) {};
\draw (1,1.5) node [draw] (l0) {};
\draw (2,1.5) node [draw] (l1) {};
\end{tikzpicture}\end{center}

We examine the cases of $d=15'$ and $d=15''$ next. As with the $H_3$-case, we simply take algebraic conjugates for the real stem, but the non-real characters are more complicated.

For $d=15''$, the tables in \cite{muller1997} and a degree argument furnish us with the real stem, but there are some twenty non-real characters remaining. The $I_2(5)$-series must be connected to the exceptional node: since all unipotent characters for $H_3(q)$ are projective (as $\Phi_{15}''$ does not divide $|H_3(q)|$) we induce $I_2(5)[1,2];\ep$ from $H_3$ to $H_4$ to get that $I_2(5)[1,2];\phi_{1,10}$ and $I_2(5)[1,2];\phi_{2,2}$ are connected, and similary for $I_2(5)[1,2];1$, and $I_2(5)[1,2];\phi_{2,2}$ and $I_2(5)[1,2];\phi_{1,0}$; a parity argument and a degree argument means that $I_2(5)[1,2];\phi_{1,10}$ is connected to the exceptional node. However, very little can be said about the locations of the cuspidal characters, except for some obvious statements with degree arguments and parity arguments. With the combinatorial Brou\'e conjecture however, we can get this diagram.
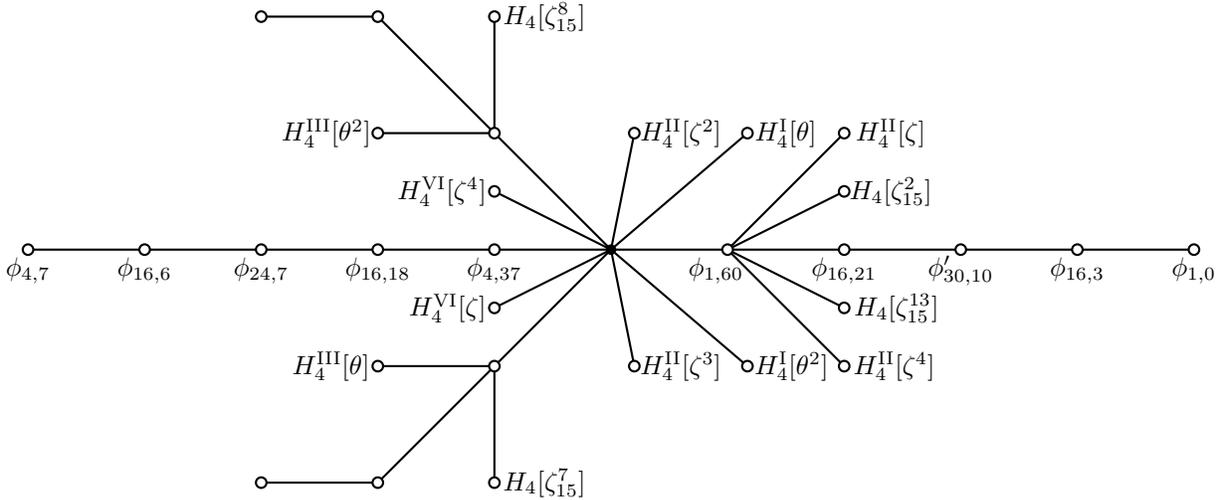
\begin{center}\begin{tikzpicture}[thick,scale=1.55]
\tikzstyle{every node}=[rectangle]

\draw (2.4,1.5) node {$H_4[\zeta_{15}^2]$};
\draw (2.4,2) node {$H_4^\mathrm{II}[\zeta]$};

\draw (2.45,0.5) node {$H_4[\zeta_{15}^{13}]$};
\draw (2.43,0) node {$H_4^\mathrm{II}[\zeta^4]$};

\draw (1.5,2) node {$H_4^\mathrm{I}[\theta]$};
\draw (1.55,0) node {$H_4^\mathrm{I}[\theta^2]$};

\draw (-0.57,3) node {$H_4[\zeta_{15}^8]$};
\draw (-0.57,-1) node {$H_4[\zeta_{15}^7]$};

\draw (-2.4,0) node {$H_4^\mathrm{III}[\theta]$};
\draw (-2.44,2) node {$H_4^\mathrm{III}[\theta^2]$};

\draw (-1.45,1.5) node {$H_4^\mathrm{VI}[\zeta^4]$};
\draw (-1.4,0.5) node {$H_4^\mathrm{VI}[\zeta]$};

\draw (0.6,2) node {$H_4^\mathrm{II}[\zeta^2]$};
\draw (0.6,0) node {$H_4^\mathrm{II}[\zeta^3]$};

\draw (5,0.82) node {$\phi_{1,0}$};
\draw (4,0.82) node {$\phi_{16,3}$};
\draw (3,0.82) node {$\phi_{30,10}'$};
\draw (2,0.82) node {$\phi_{16,21}$};
\draw (0.9,0.82) node {$\phi_{1,60}$};
\draw (-5,0.82) node {$\phi_{4,7}$};
\draw (-4,0.82) node {$\phi_{16,6}$};
\draw (-3,0.82) node {$\phi_{24,7}$};
\draw (-2,0.82) node {$\phi_{16,18}$};
\draw (-1,0.82) node {$\phi_{4,37}$};

\draw (-5,1) -- (5,1);
\draw (1,1) -- (2,0);
\draw (1,1) -- (2,2);
\draw (1,1) -- (2,0.5);
\draw (1,1) -- (2,1.5);

\draw (0,1) -- (1.17,2);
\draw (0,1) -- (1.17,0);

\draw (0,1) -- (0.2,2);
\draw (0,1) -- (0.2,0);

\draw (0,1) -- (-1,0.5);
\draw (0,1) -- (-1,1.5);

\draw (0,1) -- (-2,3);
\draw (0,1) -- (-2,-1);
\draw (-2,3) -- (-3,3);
\draw (-2,-1) -- (-3,-1);

\draw (-1,2) -- (-2,2);
\draw (-1,2) -- (-1,3);

\draw (-1,0) -- (-2,0);
\draw (-1,0) -- (-1,-1);

\tikzstyle{every node}=[circle, fill=black!0, inner sep=0pt, minimum width=4pt]
\draw (5,1) node [draw] (l0) {};
\draw (4,1) node [draw] (l0) {};
\draw (3,1) node [draw] (l1) {};
\draw (2,1) node [draw] (l1) {};
\draw (1,1) node [draw] (l1) {};
\draw (-5,1) node [draw] (l0) {};
\draw (-4,1) node [draw] (l0) {};
\draw (-3,1) node [draw] (l1) {};
\draw (-2,1) node [draw] (l1) {};
\draw (-1,1) node [draw] (l1) {};
\draw (0,1) node [fill=black!100] (ld) {};

\draw (2,1.5) node [draw] (l1) {};
\draw (2,2) node [draw] (l1) {};
\draw (2,0.5) node [draw] (l1) {};
\draw (2,0) node [draw] (l1) {};

\draw (-1,0) node [draw] (l1) {};
\draw (-1,2) node [draw] (l1) {};
\draw (-2,-1) node [draw] (l1) {};
\draw (-2,3) node [draw] (l1) {};
\draw (-3,-1) node [draw] (l1) {};
\draw (-3,3) node [draw] (l1) {};
\draw (-2,2) node [draw] (l1) {};
\draw (-1,3) node [draw] (l1) {};
\draw (-2,0) node [draw] (l1) {};
\draw (-1,-1) node [draw] (l1) {};

\draw (1.17,2) node [draw] (l1) {};
\draw (1.17,0) node [draw] (l1) {};

\draw (0.2,2) node [draw] (l1) {};
\draw (0.2,0) node [draw] (l1) {};

\draw (-1,0.5) node [draw] (l1) {};
\draw (-1,1.5) node [draw] (l1) {};
\end{tikzpicture}\end{center}
The unlabelled nodes above the real stem are, from the exceptional, $I_2(5)[1,2];\phi_{1,10}$, $I_2(5)[1,2];\phi_{2,2}$ and $I_2(5)[1,2];\phi_{1,0}$, and below are their conjugates.

For $d=15'$, a similar situation occurs. All but the cuspidal characters can be identified. Using algebraic conjugates, we get the real stem as with the previous case. Harish-Chandra induction allows us to determine the relationship between the $I_2(5)$-series characters, and a degree and parity argument proves that $I_2(5)[1,2];\phi_{1,10}$ and its dual are connected to $\phi_{1,60}$. In this case, $I_2(5)[1,2];\phi_{1,10}$ has a leading term of $q^{59}/\sqrt5$, and the only other candidate for the character to which it is connected -- $\phi_{4,31}$, as a leading term of $(5-\sqrt5)q^{59}/10$, which is smaller. For the cuspidal characters, we need to assume the combinatorial Brou\'e conjecture, and get the following.
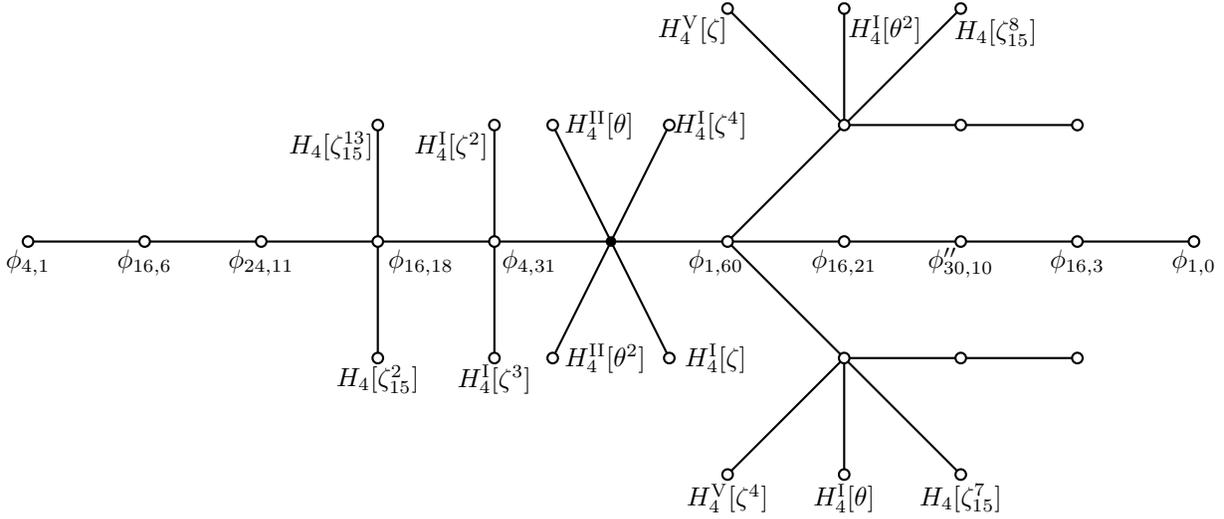
\begin{center}\begin{tikzpicture}[thick,scale=1.55]
\tikzstyle{every node}=[rectangle]

\draw (3.3,2.8) node {$H_4[\zeta_{15}^8]$};
\draw (3,-1.2) node {$H_4[\zeta_{15}^7]$};

\draw (0.85,2) node {$H_4^\mathrm{I}[\zeta^4]$};
\draw (0.9,0) node {$H_4^\mathrm{I}[\zeta]$};

\draw (-0.1,2) node {$H_4^\mathrm{II}[\theta]$};
\draw (-0.05,0) node {$H_4^\mathrm{II}[\theta^2]$};

\draw (2,-1.2) node {$H_4^\mathrm{I}[\theta]$};
\draw (1,-1.2) node {$H_4^\mathrm{V}[\zeta^4]$};

\draw (2.35,2.82) node {$H_4^\mathrm{I}[\theta^2]$};
\draw (0.7,2.82) node {$H_4^\mathrm{V}[\zeta]$};

\draw (-2,-0.18) node {$H_4[\zeta_{15}^2]$};
\draw (-2.4,1.82) node {$H_4[\zeta_{15}^{13}]$};

\draw (-1,-0.18) node {$H_4^\mathrm{I}[\zeta^3]$};
\draw (-1.37,1.82) node {$H_4^\mathrm{I}[\zeta^2]$};

\draw (5,0.82) node {$\phi_{1,0}$};
\draw (4,0.82) node {$\phi_{16,3}$};
\draw (3,0.82) node {$\phi_{30,10}''$};
\draw (2,0.82) node {$\phi_{16,21}$};
\draw (0.9,0.82) node {$\phi_{1,60}$};
\draw (-5,0.82) node {$\phi_{4,1}$};
\draw (-4,0.82) node {$\phi_{16,6}$};
\draw (-3,0.82) node {$\phi_{24,11}$};
\draw (-1.63,0.82) node {$\phi_{16,18}$};
\draw (-0.7,0.82) node {$\phi_{4,31}$};

\draw (-5,1) -- (5,1);

\draw (1,1) -- (3,3);
\draw (1,1) -- (3,-1);

\draw (2,2) -- (4,2);
\draw (2,2) -- (2,3);
\draw (2,2) -- (1,3);

\draw (2,0) -- (4,0);
\draw (2,0) -- (2,-1);
\draw (2,0) -- (1,-1);

\draw (0,1) -- (0.5,2);
\draw (0,1) -- (0.5,0);
\draw (0,1) -- (-0.5,2);
\draw (0,1) -- (-0.5,0);

\draw (-1,0) -- (-1,2);
\draw (-2,0) -- (-2,2);

\tikzstyle{every node}=[circle, fill=black!0, inner sep=0pt, minimum width=4pt]
\draw (5,1) node [draw] (l0) {};
\draw (4,1) node [draw] (l0) {};
\draw (3,1) node [draw] (l1) {};
\draw (2,1) node [draw] (l1) {};
\draw (1,1) node [draw] (l1) {};
\draw (-5,1) node [draw] (l0) {};
\draw (-4,1) node [draw] (l0) {};
\draw (-3,1) node [draw] (l1) {};
\draw (-2,1) node [draw] (l1) {};
\draw (-1,1) node [draw] (l1) {};
\draw (0,1) node [fill=black!100] (ld) {};

\draw (0.5,2) node [draw] (l1) {};
\draw (0.5,0) node [draw] (l1) {};
\draw (-0.5,2) node [draw] (l1) {};
\draw (-0.5,0) node [draw] (l1) {};

\draw (2,2) node [draw] (l1) {};
\draw (3,3) node [draw] (l1) {};
\draw (3,2) node [draw] (l1) {};
\draw (4,2) node [draw] (l1) {};
\draw (2,3) node [draw] (l1) {};
\draw (1,3) node [draw] (l1) {};

\draw (2,0) node [draw] (l1) {};
\draw (3,-1) node [draw] (l1) {};
\draw (3,0) node [draw] (l1) {};
\draw (4,0) node [draw] (l1) {};
\draw (2,-1) node [draw] (l1) {};
\draw (1,-1) node [draw] (l1) {};

\draw (-1,2) node [draw] (l1) {};
\draw (-1,0) node [draw] (l1) {};
\draw (-2,2) node [draw] (l1) {};
\draw (-2,0) node [draw] (l1) {};

%
\end{tikzpicture}\end{center}
The unlabelled nodes above the real stem are, from the exceptional, $I_2(5)[1,2];\phi_{1,10}$, $I_2(5)[1,2];\phi_{2,4}$ and $I_2(5)[1,2];\phi_{1,0}$, and below are their conjugates.

\medskip

For $20'$, $20''$, $30'$ and $30''$, the real stem consists of five principal-series characters and one other real character, toegether with the exceptional node. Induction from the $I_2(5)$-subgroup fixes the $I_2(5)$ branches, but of course not how they attach to the real stem. In each of the four cases, a simple parity and degree argument resolves the locations for the $I_2(5)$ branches, and an easy degree argument proves that the degree of the $I_2(5)$-characters increases towards the real stem. The parity and degree arguments available for cuspidal characters cannot give us much more information about their locations, and so we must assume the combinatorial Brou\'e conjecture in all four cases.

For $d=20''$ we get the following tree.
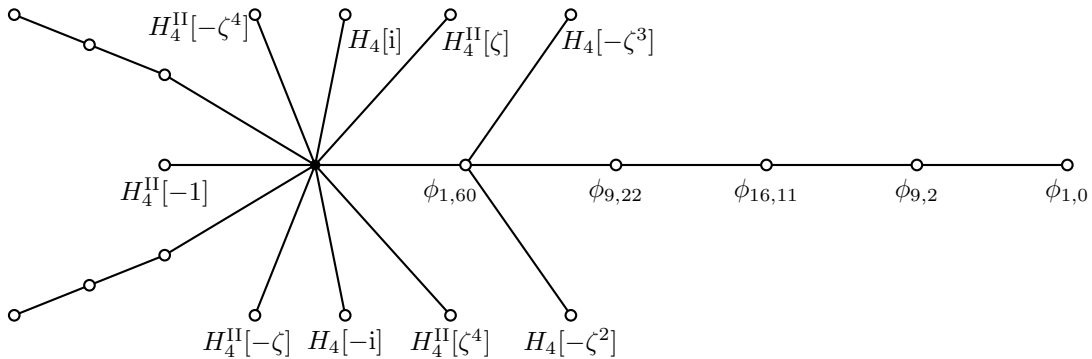
\begin{center}\begin{tikzpicture}[thick,scale=2]

\tikzstyle{every node}=[rectangle]

\draw (5,0.82) node {$\phi_{1,0}$};
\draw (4,0.82) node {$\phi_{9,2}$};
\draw (3,0.82) node {$\phi_{16,11}$};
\draw (2,0.82) node {$\phi_{9,22}$};
\draw (0.9,0.82) node {$\phi_{1,60}$};
\draw (-1,0.82) node {$H_4^\mathrm{II}[-1]$};

\draw (0.4,1.82) node {$H_4[\I]$};
\draw (0.2,-0.18) node {$H_4[-\I]$};

\draw (1.07,1.82) node {$H_4^\mathrm{II}[\zeta]$};
\draw (0.9,-0.18) node {$H_4^\mathrm{II}[\zeta^4]$};

\draw (1.95,1.82) node {$H_4[-\zeta^3]$};
\draw (1.7,-0.18) node {$H_4[-\zeta^2]$};

\draw (-0.45,-0.18) node {$H_4^\mathrm{II}[-\zeta]$};
\draw (-0.77,1.92) node {$H_4^\mathrm{II}[-\zeta^4]$};

\tikzstyle{every node}=[circle, fill=black!0, inner sep=0pt, minimum width=4pt]

\draw (-1,1) -- (5,1);
\draw (0,1) -- (0.2,2);
\draw (0,1) -- (0.2,0);

\draw (0,1) -- (0.9,2);
\draw (0,1) -- (0.9,0);
\draw (1,1) -- (1.7,2);
\draw (1,1) -- (1.7,0);

\draw (0,1) -- (-0.4,2);
\draw (0,1) -- (-0.4,0);

\draw (0,1) -- (-1,0.4);
\draw (0,1) -- (-1,1.6);

\draw (-1,0.4) -- (-2,0);
\draw (-1,1.6) -- (-2,2);

\draw (5,1) node [draw] (l0) {};
\draw (4,1) node [draw] (l0) {};
\draw (3,1) node [draw] (l1) {};
\draw (2,1) node [draw] (l1) {};
\draw (1,1) node [draw] (l1) {};
\draw (-1,1) node [draw] (l1) {};
\draw (0,1) node [fill=black!100] (ld) {};

\draw (0.2,0) node [draw] (l0) {};
\draw (0.2,2) node [draw] (l0) {};
\draw (0.9,0) node [draw] (l0) {};
\draw (0.9,2) node [draw] (l0) {};
\draw (1.7,0) node [draw] (l0) {};
\draw (1.7,2) node [draw] (l0) {};

\draw (-0.4,0) node [draw] (l0) {};
\draw (-0.4,2) node [draw] (l0) {};

\draw (-1,1.6) node [draw] (l0) {};
\draw (-1,0.4) node [draw] (l0) {};
\draw (-1.5,1.8) node [draw] (l0) {};
\draw (-1.5,0.2) node [draw] (l0) {};
\draw (-2,0) node [draw] (l0) {};
\draw (-2,2) node [draw] (l0) {};
\end{tikzpicture}\end{center}
The unlabelled nodes above the real stem are, from the exceptional, $I_2(5)[1,3];\phi_{1,5}'$, $I_2(5)[1,3];\phi_{2,4}$ and $I_2(5)[1,3];\phi_{1,5}''$, and their conjugates are below. This tree follows another pattern shown by the exceptional groups of Lie type in the `sub-Coxeter case', i.e., the $\Phi_d$ such that $d$ is second-largest: in these cases, all non-real characters are on branches emanating from the exceptional node  -- arranged in increasing argument of the eigenvalue of the Frobenius, as for the Coxeter case -- apart from one pair of cuspidal characters (should any exist), which emanate from the Steinberg character. Notice that we also see this with $H_3(q)$ and $d=6$. Potentially there is a `sub-Coxeter theory' to go along with the theory associated with the Coxeter torus for any group of Lie type.

For $d=20'$ however, the tree is more complicated, indicating a general trend for the single-primed cyclotomic polynomials to have more complex trees than the double-primed ones.
\begin{center}\begin{tikzpicture}[thick,scale=2]
\tikzstyle{every node}=[rectangle, fill=black!0, inner sep=0pt, minimum width=4pt]
\draw (5,0.82) node {$\phi_{1,0}$};
\draw (4,0.82) node {$\phi_{9,6}$};
\draw (3,0.82) node {$\phi_{16,13}$};
\draw (2,0.82) node {$\phi_{9,26}$};
\draw (1.22,0.74) node {$\phi_{1,60}$};
\draw (-1,0.82) node {$H_4^\mathrm{III}[-1]$};

\draw (0,0.32) node {$H_4^\mathrm{I}[\zeta]$};
\draw (1.34,2) node {$H_4[-\I]$};
\draw (-0.3,1.82) node {$H_4^\mathrm{I}[-\zeta^2]$};
\draw (2,0.32) node {$H_4^\mathrm{I}[-\zeta]$};
\draw (3.7,1.97) node {$I_2(5)[1,2];\phi_{1,5}''$};
\draw (2.5,1.82) node {$I_2(5)[1,2];\phi_{2,2}'$};

\draw (-0.2,1.32) node {$H_4^\mathrm{I}[\zeta^4]$};
\draw (0,-0.18) node {$H_4^\mathrm{I}[-\zeta^3]$};
\draw (1,-0.18) node {$H_4[\I]$};
\draw (1.9,-0.18) node {$I_2(5)[1,3];\phi_{2,2}'$};
\draw (3.2,-0.18) node {$I_2(5)[1,3];\phi_{1,5}''$};
\draw (2.23,1.32) node {$H_4^\mathrm{I}[-\zeta^4]$};
\tikzstyle{every node}=[circle, fill=black!0, inner sep=0pt, minimum width=4pt]
\draw (-1,1) -- (5,1); 
\draw (1,0) -- (1,2); 

\draw (1,1) -- (2,1.5);
\draw (1,1) -- (2,0.5);

\draw (1,1) -- (0,1.5);
\draw (1,1) -- (0,0.5);

\draw (1,1.5) -- (2,2);
\draw (2,2) -- (3,2);
\draw (2,0) -- (3,0);
\draw (1,1.5) -- (0,2);
\draw (1,0.5) -- (2,0);
\draw (1,0.5) -- (0,0);

\draw (5,1) node [draw] (l0) {};
\draw (4,1) node [draw] (l0) {};
\draw (3,1) node [draw] (l1) {};
\draw (2,1) node [draw] (l1) {};
\draw (1,1) node [draw] (l1) {};
\draw (-1,1) node [draw] (l1) {};
\draw (0,1) node [fill=black!100] (ld) {};
\draw (0,0.5) node [draw] (l0) {};
\draw (1,0.5) node [draw] (l0) {};
\draw (2,0.5) node [draw] (l1) {};
\draw (0,1.5) node [draw] (l0) {};
\draw (1,1.5) node [draw] (l0) {};
\draw (2,1.5) node [draw] (l1) {};
\draw (0,2) node [draw] (l0) {};
\draw (1,2) node [draw] (l0) {};
\draw (2,2) node [draw] (l1) {};
\draw (3,2) node [draw] (l1) {};
\draw (0,0) node [draw] (l0) {};
\draw (1,0) node [draw] (l0) {};
\draw (2,0) node [draw] (l1) {};
\draw (3,0) node [draw] (l1) {};
\end{tikzpicture}\end{center}
The unlabelled node above the real stem is $I_2(5)[1,2];\phi_{1,5}'$, and its conjugate labels the conjugate node.

\medskip

For $d=30'$ we get the most complicated tree for any (genuine or not) group of Lie type.
\begin{center}\begin{tikzpicture}[thick,scale=2]

\draw (5,0.82) node {$\phi_{1,0}$};
\draw (4,0.82) node {$\phi_{4,7}$};
\draw (3,0.82) node {$\phi_{6,20}$};
\draw (2,0.82) node {$\phi_{4,37}$};
\draw (1.22,0.74) node {$\phi_{1,60}$};
\draw (-1,0.82) node {$H_4^\mathrm{IV}[-1]$};
\draw (1.62,2.5) node {$I_2(5)[1,3];\phi_{1,0}$};
\draw (0.4,2.5) node {$H_3[-\I];1$};
\draw (-0.67,2.5) node {$H_4[\zeta_{15}^{13}]$};
\draw (-2.25,2.32) node {$H_4^\mathrm{II}[-\zeta^3]$};
\draw (-2,1.82) node {$H_4^\mathrm{II}[\theta]$};
\draw (-1.25,1.82) node {$H_3[-\I];\ep$};
\draw (1.37,2) node {$H_4[-\theta^2]$};
\draw (2.23,1.32) node {$H_4[\zeta_{15}^7]$};
\draw (1.38,1.5) node {$H_4^\mathrm{II}[-\zeta]$};
\draw (-1,1.32) node {$H_4^\mathrm{III}[\zeta^4]$};
\draw (1,0.32) node {$H_4^\mathrm{II}[-\zeta^4]$};
\draw (2,0.32) node {$H_4[\zeta_{15}^8]$};
\draw (-1,0.32) node {$H_4^\mathrm{III}[\zeta]$};
\draw (-2,-0.18) node {$H_4^\mathrm{II}[\theta^2]$};
\draw (-0.55,0) node {$H_3[\I];\ep$};
\draw (1,-0.68) node {$I_2(5)[1,2];\phi_{1,0}$};
\draw (1,-0.18) node {$H_4[-\theta]$};
\draw (0,-0.68) node {$H_3[\I];1$};
\draw (-1,-0.68) node {$H_4[\zeta_{15}^2]$};
\draw (-2,-0.68) node {$H_4^\mathrm{II}[-\zeta^2]$};

\draw (-1,1) -- (5,1); 
\draw (1,0.5) -- (1,1.5); 

\draw (2,0.5) -- (1,1);
\draw (1,1) -- (2,1.5);
\draw (1,1) -- (-2,2.5);
\draw (1,1) -- (-2,-0.5);

\draw (0,1.5) -- (-1,1.5);
\draw (-1,2) -- (-2,2);

\draw (0,1.5) -- (0,2);
\draw (0,1.5) -- (1,2);
\draw (-1,2) -- (-1,2.5);
\draw (-1,2) -- (0,2.5);
\draw (0,2) -- (1,2.5);

\draw (0,0.5) -- (-1,0.5);
\draw (-1,0) -- (-2,0);

\draw (0,0.5) -- (0,0);
\draw (0,0.5) -- (1,0);
\draw (-1,0) -- (-1,-0.5);
\draw (-1,0) -- (0,-0.5);
\draw (0,0) -- (1,-0.5);

\draw (5,1) node [draw] (l0) {};
\draw (4,1) node [draw] (l0) {};
\draw (3,1) node [draw] (l1) {};
\draw (2,1) node [draw] (l1) {};
\draw (1,1) node [draw] (l1) {};
\draw (-1,1) node [draw] (l1) {};
\draw (0,1) node [fill=black!100] (ld) {};

\draw (2,1.5) node [draw] (l0) {};
\draw (1,1.5) node [draw] (l1) {};
\draw (0,1.5) node [draw] (l1) {};
\draw (-1,2) node [draw] (l1) {};
\draw (1,2) node [draw] (l0) {};
\draw (0,2) node [draw] (l1) {};
\draw (-1,1.5) node [draw] (l1) {};
\draw (-2,2) node [draw] (l1) {};
\draw (1,2.5) node [draw] (l0) {};
\draw (0,2.5) node [draw] (l1) {};
\draw (-1,2.5) node [draw] (l1) {};
\draw (-2,2.5) node [draw] (l1) {};

\draw (2,0.5) node [draw] (l0) {};
\draw (1,0.5) node [draw] (l1) {};
\draw (0,0.5) node [draw] (l1) {};
\draw (-1,0.5) node [draw] (l1) {};
\draw (1,0) node [draw] (l0) {};
\draw (0,0) node [draw] (l1) {};
\draw (-1,0) node [draw] (l1) {};
\draw (-2,0) node [draw] (l1) {};
\draw (1,-0.5) node [draw] (l0) {};
\draw (0,-0.5) node [draw] (l1) {};
\draw (-1,-0.5) node [draw] (l1) {};
\draw (-2,-0.5) node [draw] (l1) {};
\end{tikzpicture}\end{center}
The unlabelled nodes above the real stem are (in order) $I_2(5)[1,3];\phi_{1,10}$ and $I_2(5)[1,3];\phi_{2,3}$, and their conjugates label the conjugate nodes.

\medskip

Finally, for $d=30''$ the tree is consistent with the Hiss--L\"ubeck--Malle conjecture, and we do not draw the tree explicitly but refer to the description given in Section \ref{sec:Brauertreeargs}.

\medskip

We now turn to non-principal blocks, as with $H_3$. Again, there are two non-principal blocks with cyclic defect groups for $d=1$, this time with $1$-cuspidal pairs $(H_3.\Phi_1,H_3[\pm\I])$, and consisting of a line with exceptional node in the middle, and $H_3[\pm\I];1$ on the one side and $H_3[\pm\I];\ep$ on the other. For $d=2$ there are again two non-principal blocks with cyclic defect groups, this time with $d$-split Levi $H_3.\Phi_2$ and unipotent characters $\phi_{4,3}$ and $\phi_{4,4}$. The Brauer trees are lines with the exceptional node at the end, with $\phi_{16,18}$ connected to the exceptional and $\phi_{16,3}$ connected to that for the character $\phi_{4,3}$, and $\phi_{16,21}$ connected to the exceptional and $\phi_{16,6}$ connected to that for the character $\phi_{4,4}$.

Next we consider the case $d=5'$. Although the Brauer tree is the same as the corresponding $d$ for $H_3$, Harish-Chandra induction is \emph{not} a Morita equivalence. However, the Harish-Chandra induction of $\phi_{1,15}+I_2(5)[1,3];\ep$ or $\phi_{3,6}+I_2(5)[1,3];\ep$ -- one of which we know to be projective even without the combinatorial Brou\'e conjecture, and with no restriction on $q$ -- when cut by this particular non-principal block, becomes $\phi_{4,31}+\phi_{9,22}+I_2(5)[1,3];\phi_{2,1}+I_2(5)[1,3];\phi_{2,4}$. Since this character is projective, this means that exactly one of $I_2(5)[1,3];\phi_{2,1}$ and $I_2(5)[1,3];\phi_{2,4}$ is connected to $\phi_{4,31}$, and the other is connected to $\phi_{9,22}$. We cannot label the non-real nodes of the Brauer tree using this process, never mind get the full planar embedding, but we \emph{can} determine the shape of the Brauer tree.
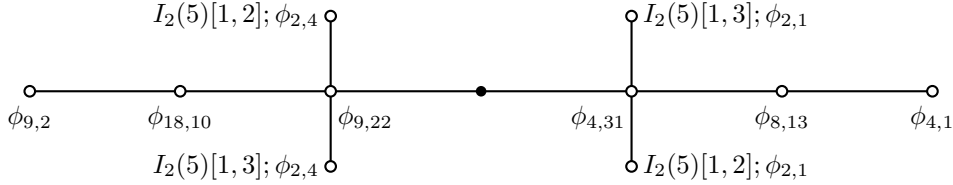
\begin{center}\begin{tikzpicture}[thick,scale=2]

\tikzstyle{every node}=[rectangle]

\draw (4.63,0.5) node {$I_2(5)[1,2];\phi_{2,1}$};
\draw (1.37,1.5) node {$I_2(5)[1,2];\phi_{2,4}$};
\draw (4.63,1.5) node {$I_2(5)[1,3];\phi_{2,1}$};
\draw (1.37,0.5) node {$I_2(5)[1,3];\phi_{2,4}$};

\draw (6,0.82) node {$\phi_{4,1}$};
\draw (5,0.82) node {$\phi_{8,13}$};
\draw (3.78,0.82) node {$\phi_{4,31}$};
\draw (2.23,0.82) node {$\phi_{9,22}$};
\draw (1,0.82) node {$\phi_{18,10}$};
\draw (0,0.82) node {$\phi_{9,2}$};

\tikzstyle{every node}=[circle, fill=black!0, inner sep=0pt, minimum width=4pt]

\draw (0,1) -- (6,1);
\draw (4,0.5) -- (4,1.5);
\draw (2,0.5) -- (2,1.5);

\draw (6,1) node [draw] (l0) {};
\draw (5,1) node [draw] (l1) {};
\draw (4,1) node [draw] (l1) {};
\draw (0,1) node [draw] (l1) {};
\draw (3,1) node [fill=black!100] (ld) {};
\draw (2,1) node [draw] (l0) {};
\draw (1,1) node [draw] (l1) {};
\draw (2,0.5) node [draw] (l0) {};
\draw (2,1.5) node [draw] (l1) {};
\draw (4,0.5) node [draw] (l0) {};
\draw (4,1.5) node [draw] (l1) {};
\end{tikzpicture}\end{center}

As mentioned in the previous section, this new tree can be used to pin down the structure of the Brauer tree for the principal block of $H_3(q)$ when $d=5'$ as well. We now take Harish-Chandra restriction of $\phi_{3,4}+I_2(5)[1,3];\phi_{2,1}$ or $\phi_{3,4}+I_2(5)[1,3];\phi_{2,4}$ -- one of which is projective regardless of the labelling (which we couldn't determine without the combinatorial Brou\'e conjecture) -- to $H_3$ and cut by the principal block, leaving $\phi_{1,15}+\phi_{3,6}+I_2(5)[1,3];1+I_2(5)[1,3];\ep$ in both cases. This means that the two $I_2(5)$-branches cannot both be connected to $\phi_{1,15}$, so one is connected to $\phi_{3,6}$.

\medskip

For $d=5''$, again the block is Morita equivalent to the principal block of $H_3(q)$ for $d=5''$, but again Harish-Chandra induction does not induce a Morita equivalence. The real stem consists solely of principal series characters, so a M\"uller argument and a degree argument gives us six of the ten unipotent characters. A parity argument means that the non-real characters -- each with leading term in the degree $(1+\sqrt5)q^{54}/20$ are connected either to the exceptional, or one $\phi_{8,13}$ and $\phi_{18,10}$, each of which has leading term in the degree $q^{54}/10$; this proves that all non-real characters are connected to the exceptional, but as with the $H_3$ case does not compute the planar embedding. Assuming the combinatorial Brou\'e conjecture we get the following.
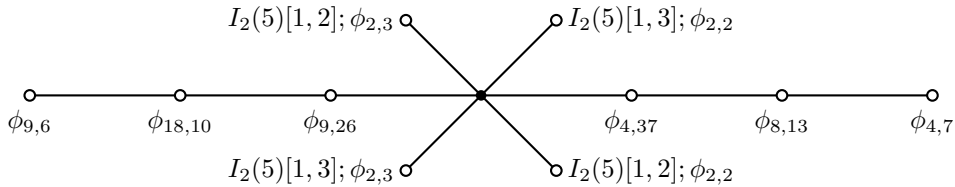
\begin{center}\begin{tikzpicture}[thick,scale=2]
\tikzstyle{every node}=[rectangle]
\draw (4.13,0.5) node {$I_2(5)[1,2];\phi_{2,2}$};
\draw (1.87,1.5) node {$I_2(5)[1,2];\phi_{2,3}$};
\draw (4.13,1.5) node {$I_2(5)[1,3];\phi_{2,2}$};
\draw (1.87,0.5) node {$I_2(5)[1,3];\phi_{2,3}$};

\draw (6,0.82) node {$\phi_{4,7}$};
\draw (5,0.82) node {$\phi_{8,13}$};
\draw (4,0.82) node {$\phi_{4,37}$};
\draw (2,0.82) node {$\phi_{9,26}$};
\draw (1,0.82) node {$\phi_{18,10}$};
\draw (0,0.82) node {$\phi_{9,6}$};

\tikzstyle{every node}=[circle, fill=black!0, inner sep=0pt, minimum width=4pt]

\draw (0,1) -- (6,1);
\draw (3.5,0.5) -- (2.5,1.5);
\draw (2.5,0.5) -- (3.5,1.5);

\draw (6,1) node [draw] (l0) {};
\draw (5,1) node [draw] (l1) {};
\draw (4,1) node [draw] (l1) {};
\draw (0,1) node [draw] (l1) {};
\draw (3,1) node [fill=black!100] (ld) {};
\draw (2,1) node [draw] (l0) {};
\draw (1,1) node [draw] (l1) {};
\draw (2.5,0.5) node [draw] (l0) {};
\draw (2.5,1.5) node [draw] (l1) {};
\draw (3.5,0.5) node [draw] (l0) {};
\draw (3.5,1.5) node [draw] (l1) {};
\end{tikzpicture}\end{center}

We now use this to prove that $I_2(5)[1,3];1$ and its dual are connected to the exceptional node for $H_3(q)$ and $d=5''$: if they were not then by a parity argument they are connected to either $\phi_{4,3}$ or $\phi_{4,4}$: Harish-Chandra induction of either $\phi_{4,3}+I_2(5)[1,3];1$ or $\phi_{4,4}+I_2(5)[1,3];1$ to $H_4(q)$, cut by the non-principal block above, yields $\phi_{8,13}+\phi_{18,1-}+I_2(5)[1,3];\phi_{2,2}+I_2(5)[1,3];\phi_{2,3}$, which is definitely not projective. This completes the proof for $d=5''$.

\medskip

Finally we have $d=10'$ and $d=10''$. Six of the unipotent characters lie on the real stem in each case, and so a M\"uller argument and a degree argument give us the real stem. For $d=10'$, by a parity argument the four non-real cuspidal characters, which have leading term of their degrees $(1+\sqrt5)q^{54}/20$. must lie on either $\phi_{4,37}$, or either $H_4^\mathrm{I}[-1]$ or $\phi_{10,12}$, the latter two having leading term of the degre $q^{54}/10$. This gives the shape of the Brauer tree below, but as usual the planar embedding requires the combinatorial Brou\'e conjecture. (If $q=2$ then the cuspidals could also be attached to $\phi_{10,12}$.)

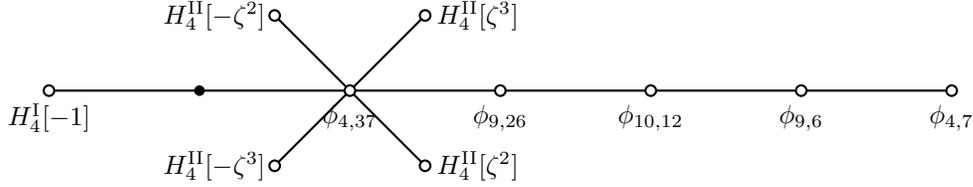
\begin{center}\begin{tikzpicture}[thick,scale=2]
\tikzstyle{every node}=[rectangle]
\draw (2.85,0.5) node {$H_4^\mathrm{II}[\zeta^2]$};
\draw (1.1,1.5) node {$H_4^\mathrm{II}[-\zeta^2]$};
\draw (2.85,1.5) node {$H_4^\mathrm{II}[\zeta^3]$};
\draw (1.1,0.5) node {$H_4^\mathrm{II}[-\zeta^3]$};

\draw (6,0.82) node {$\phi_{4,7}$};
\draw (5,0.82) node {$\phi_{9,6}$};
\draw (4,0.82) node {$\phi_{10,12}$};
\draw (3,0.82) node {$\phi_{9,26}$};
\draw (2,0.82) node {$\phi_{4,37}$};
\draw (0,0.82) node {${H_4^\mathrm{I}[-1]}$};

\tikzstyle{every node}=[circle, fill=black!0, inner sep=0pt, minimum width=4pt]

\draw (0,1) -- (6,1);
\draw (2.5,0.5) -- (1.5,1.5);
\draw (1.5,0.5) -- (2.5,1.5);

\draw (6,1) node [draw] (l0) {};
\draw (5,1) node [draw] (l1) {};
\draw (4,1) node [draw] (l1) {};
\draw (0,1) node [draw] (l1) {};
\draw (1,1) node [fill=black!100] (ld) {};
\draw (3,1) node [draw] (l0) {};
\draw (2,1) node [draw] (l1) {};
\draw (2.5,0.5) node [draw] (l0) {};
\draw (2.5,1.5) node [draw] (l1) {};
\draw (1.5,0.5) node [draw] (l0) {};
\draw (1.5,1.5) node [draw] (l1) {};
\end{tikzpicture}\end{center}

When $d=10''$ however, there is more than one acceptable location for the non-real characters, both the exceptional and $\phi_{9,22}$. Indeed, according to the combinatorial Brou\'e conjecture, one pair of non-real characters emanates from each of the potential nodes on the real stem that can be adjacent to non-real characters, giving the tree below.
\begin{center}\begin{tikzpicture}[thick,scale=2]
\tikzstyle{every node}=[rectangle]

\draw (2.4,1.5) node {$H_4^\mathrm{I}[-\zeta^3]$};
\draw (2.4,0.5) node {$H_4^\mathrm{I}[-\zeta^2]$};
\draw (0.3,1.5) node {$H_4^\mathrm{I}[\zeta^2]$};
\draw (0.3,0.5) node {$H_4^\mathrm{I}[\zeta^3]$};

\draw (5,0.82) node {$\phi_{4,1}$};
\draw (4,0.82) node {$\phi_{9,2}$};
\draw (3,0.82) node {$\phi_{10,12}$};
\draw (2.22,0.82) node {$\phi_{9,22}$};
\draw (1,0.82) node {$\phi_{4,31}$};
\draw (-1,0.82) node {$H_4^\mathrm{I}[-1]$};

\tikzstyle{every node}=[circle, fill=black!0, inner sep=0pt, minimum width=4pt]

\draw (-1,1) -- (5,1);
\draw (2,0.5) -- (2,1.5);
\draw (0,0.5) -- (0,1.5);

\draw (5,1) node [draw] (l0) {};
\draw (4,1) node [draw] (l0) {};
\draw (3,1) node [draw] (l1) {};
\draw (2,1) node [draw] (l1) {};
\draw (1,1) node [draw] (l1) {};
\draw (0,1) node [fill=black!100] (ld) {};
\draw (-1,1) node [draw] (l1) {};
\draw (0,0.5) node [draw] (l0) {};
\draw (0,1.5) node [draw] (l1) {};

\draw (2,0.5) node [draw] (l0) {};
\draw (2,1.5) node [draw] (l1) {};
\end{tikzpicture}\end{center}

\providecommand{\bysame}{\leavevmode\hbox to3em{\hrulefill}\thinspace}

\end{document}